\documentclass[12pt]{amsart}

\usepackage{amsmath, amsthm, amssymb}
\usepackage[top=1.25in, bottom=1.25in, left=1.0in, right=1.0in]{geometry}

\usepackage{verbatim}
\usepackage{url}

\usepackage{tikz,tkz-graph}
\usepackage{subfig}
\usepackage{hyperref}

\theoremstyle{plain}
\newtheorem{thm}{Theorem}

\newtheorem*{maintheorem}{Main Theorem}

\newtheorem{fact}{Fact}
\newtheorem{lem}[thm]{Lemma}

\newtheorem*{conj}{Conjecture}

\theoremstyle{definition}

\theoremstyle{remark}

\newcommand{\chif}{\chi_{f}}
\newcommand{\ch}{{\textrm{ch}}}

\newcommand{\floor}[1]{\left\lfloor#1\right\rfloor}

\newcommand{\set}[1]{\left\{ #1 \right\}}

\author{Daniel W. Cranston}\thanks{Virginia Commonwealth University, Richmond,
VA.  \texttt{dcranston@vcu.edu}} 
\author{Landon Rabern}\thanks{LBD Data Solutions, Lancaster, PA.
\texttt{landon.rabern@gmail.com}}
\begin{document}
\title{\mbox{Planar graphs are $9/2$-colorable} 
%\mbox{and have independence ratio at least $3/13$}
}
\begin{abstract}
We show that every planar graph $G$ has a 2-fold 9-coloring.  In particular,
this implies that $G$ has fractional chromatic number at most $\frac92$.
This is the first proof (independent of the 4~Color Theorem) that there exists a
constant $k<5$ such that every planar $G$ has fractional chromatic number at
most $k$.  
%We also show that every $n$-vertex planar graph has an independent set of size
%at least $\frac{3n}{13}$.  This improves on Albertson's bound of $\frac{2n}9$.
\end{abstract}
\maketitle

\section{Introduction}

To fractionally color a graph $G$, we assign to each independent set in
$G$ a nonnegative weight, such that for each vertex $v$ the sum of the weights
on the independent sets containing $v$ is 1.  A graph $G$ is
\emph{fractionally $k$-colorable} if $G$ has such an assignment of
weights where the sum of the weights is at most $k$.  The minimum
$k$ such that $G$ is fractionally $k$-colorable is its \emph{fractional
chromatic number}, denoted $\chif(G)$.
(If we restrict the weight on each independent set to be either 0 or 1, we return
to the standard definition of chromatic number.)
In 1997, Scheinerman and Ullman~\cite[p.~75]{SU-book}
succintly described the state of the art for fractionally coloring planar graphs.  Not
much has changed since then.
\begin{quotation}
The fractional analogue of the four-color theorem is the assertion that the
maximum value of $\chif(G)$ over all planar graphs $G$ is 4. That this maximum
is no more than 4 follows from the four-color theorem itself, while the example
of $K_4$ shows that it is no less than 4. Given that the proof of the four-color
theorem is so difficult, one might ask whether it is possible to prove an
interesting upper bound for this maximum without appeal to the four-color
theorem. Certainly $\chif(G) \le 5$ for any planar $G$, because $\chi(G) \le
5$, a result whose proof is elementary.  But what about a simple proof of, say,
$\chif(G) \le \frac92$ for all planar G? The only result in this direction is
in a 1973 paper of Hilton, Rado, and Scott~\cite{HRS} that predates the proof
of the four-color theorem; they prove $\chif(G) < 5$ for any planar graph G,
although they are not able to find any constant $c < 5$ with $\chif(G) < c$ for
all planar graphs G. This may be the first appearance in print of the invariant
$\chif$.
\end{quotation}

In Section~\ref{fractional}, we give exactly what Scheinerman and Ullman asked
for---a simple proof that $\chif(G)\le \frac92$ for every planar graph $G$.  In
fact, this result is an immediate corollary of a stronger statement in our
main theorem.  Before we can express it precisely, we need another definition.
A \emph{$k$-fold $\ell$-coloring} of a graph $G$ assigns to each vertex a set of $k$
colors, such that adjacent vertices receive disjoint sets, and the union of all
sets has size at most $\ell$.  If $G$ has a $k$-fold $\ell$-coloring, then
$\chif(G)\le \frac{\ell}{k}$.  To see this, consider the $\ell$ independent sets
induced by the color classes; assign to each of these sets the weight
$\frac1{k}$.  Now we can state the theorem.

\begin{maintheorem}
Every planar graph $G$ has a 2-fold 9-coloring.  In particular,
$\chif(G)\le\frac92$.
\label{mainthm}
\end{maintheorem}

In an intuitive sense, the Main Theorem sits somewhere between the 4
Color Theorem and the 5 Color Theorem.  It is certainly implied by the
former, but it does not immediately imply the latter.  
The \emph{Kneser graph $K_{n:k}$} has as its vertices the $k$-element subsets of
$\set{1,\ldots,n}$ and two vertices are adjacent if their corresponding sets are
disjoint.  Saying that a graph $G$ has a 2-fold 9-coloring is equivalent to
saying that it has a homomorphism to the Kneser graph $K_{9:2}$.  To claim that
a coloring result for planar graphs is between the 4 and 5 Color Theorems,
we would like to show that every planar graph $G$ has a homomorphism to a graph
$H$, such that $H$ has clique number 4 and chromatic number 5.  Unfortunately,
$K_{9:2}$ is not such a graph.  It is easy to see that
$\omega(K_{n:k})=\floor{n/k}$; so $\omega(K_{9:2})=4$, as desired.  However,
Lov\'asz~\cite{Lovasz78} showed that $\chi(K_{n:k})=n-2k+2$; thus
$\chi(K_{9:2})=9-2(2)+2=7$.  Fortunately, we can easily overcome this problem.

The \emph{categorical product} (or \emph{universal product}) of graphs $G_1$ and
$G_2$, denoted $G_1\times G_2$ is defined as follows.  Let $V(G_1\times G_2) =
\{(u,v) | u \in V(G_1)\mbox{ and } v\in V(G_2)\}$; now $(u_1,v_1)$ is adjacent
to $(u_2,v_2)$ if $u_1u_2\in E(G_1)$ and $v_1v_2\in E(G_2)$.
Let $H=K_5\times K_{9:2}$.  It is well-known~\cite{Hell79} that if a graph $G$
has a homomorphism to each of graphs $G_1$ and $G_2$, then $G$ also has a
homomorphism to $G_1\times G_2$ (the image of each vertex in the product is
just the products of its images in $G_1$ and $G_2$).  The 5 Color Theorem says
that every planar graph has a homomorphism to $K_5$; so if we prove that every
planar graph $G$ has a homomorphism to $K_{9:2}$, then we also get that $G$ has
a homomorphism to $K_5\times K_{9:2}$.

It is easy to check that for any $G_1$
and $G_2$, we have $\omega(G_1\times G_2)=\min(\omega(G_1),\omega(G_2))$ and
$\chi(G_1\times G_2)\le \min(\chi(G_1),\chi(G_2))$.  
To prove this inequality, we simply color each vertex $(u,v)$ of the product with the
color of $u$ in an optimal coloring of $G_1$, or the color of $v$ in an optimal
coloring of $G_2$.  (It is an open problem whether this inequality always holds
with equality~\cite{Sauer01}.) When $H=K_5\times K_{9:2}$ we get $\omega(H)=4$ and
$\chi(H)=5$.  Earlier work of Naserasr~\cite{Naserasr06} and Ne\v{s}et\v{r}il and Ossona
de Mendez~\cite{NOdM06} also constructed graphs $H$, with $\omega(H)=4$ and
$\chi(H)=5$,  such that every planar graph $G$ has a homomorphism to $H$;
however, their examples had more vertices than ours.  Naserasr gave a graph with
size $63{62 \choose 4} = 35,144,235$ and the construction in~\cite{NOdM06}
was still larger.  In contrast, $|K_5\times K_{9:2}| = 5{9\choose 2} = 180$.

Wagner~\cite{Wagner37} characterized $K_5$-minor-free graphs. 
The Wagner graph is formed from an 8-cycle by adding an edge joining each pair
of vertices that are distance 4 along the cycle.  Wagner showed that every
maximal $K_5$-minor-free graph can be formed recursively from planar graphs and
copies of the Wagner graph by pasting along copies of $K_2$ and $K_3$  (see
also~\cite[p.~175]{Diestel}).  Since
the Wagner graph is 3-colorable, it clearly has a 2-fold 9-coloring.  To show
that every $K_5$-minor-free graph is 2-fold 9-colorable, we color each smaller
planar graph and copy of the Wagner graph, then permute colors so that the
colorings agree on the vertices that are pasted together.

Haj\'{o}s conjectured that every graph is $(k-1)$-colorable unless it contains a
subdivision of $K_k$.  This is known to be true for $k\le 4$ and false for $k\ge
7$.  The cases $k=5$ and $k=6$ remain unresolved.  Since this problem seems
difficult, we offer the following weaker conjecture. 

\begin{conj}
Every graph with no $K_5$-subdivision is 2-fold 9-colorable.
\end{conj}

%In Section~\ref{independence} we study the independence number of planar graphs.
An immediate consequence of the 4 Color Theorem is that every
$n$-vertex planar graph has an independent set of size at least $\frac{n}4$ (and this
is best possible, as shown by the disjoint union of many copies of $K_4$).
In 1968, Erd\"os~\cite{berge} suggested that perhaps this corollary could be
proved more easily than the full 4 Color Theorem.  And in 1976,
Albertson~\cite{Albertson} showed (independently of the 4~Color Theorem)
that every $n$-vertex planar graph has an independent set of size at least
$\frac{2n}9$.  

%Our second theorem improves this bound to $\frac{3n}{13}$.
%
%\begin{thm}
%Every $n$-vertex planar graph has an independent set of size at least
%$\frac{3n}{13}$.
%\label{mainthm2}
%\end{thm}

%The proof of Theorem~\ref{mainthm2} is heavily influenced by Albertson's proof.
%One apparent difference is that our proof uses the discharging method, while his
%does not.  However, this distinction is largely cosmetic.  To demonstrate this
%point, we include an appendix with a short discharging version of the final
%step in Albertson's proof, which he verified using edge-counting (the reader
%unfamiliar with discharging arguments may prefer to start with this appendix). 
%Although the arguments are essentially equivalent, the discharging method is
%somewhat more flexible.  In part it was this added flexibility that allowed us
%to push his ideas further.

Albertson's proof inspired and heavily influenced our proof of the Main Theorem.
%The proof of our Main Theorem has the following general outline.  
The bulk of the work in our proof consists in showing that certain configurations
are \emph{reducible}, i.e., they cannot appear in a minimal counterexample to
the theorem.  The proof concludes via a discharging argument,
where we show that every planar graph contains one of the forbidden
configurations; hence, it is not a minimal counterexample.  
%In the proof of the fractional coloring bound, the details are much simpler, so
%we present that result first.

Before the proof, we need a few definitions.
A \emph{$k$-vertex} is a vertex of degree $k$; similarly, a $k^-$-vertex
(resp.~$k^+$-vertex) has degree at most (resp.~at least) $k$.
A $k$-neighbor of a vertex $v$ is a $k$-vertex that is a neighbor of $v$; and
$k^-$-neighbors and $k^+$-neigbors are defined analogously.
A $k$-cycle is a cycle of length $k$.
A vertex set $V_1$ in a connected graph $G$ is \emph{separating} if $G\setminus V_1$
has at least two components.  A cycle $C$ is separating if $V(C)$ is separating.
Finally, an \emph{independent $k$-set} is an
independent set (or stable set) of size $k$.

\section{Fractional Coloring of Planar Graphs}
\label{fractional}

Now we prove our Main Theorem, that every planar graph has a 2-fold 9-coloring.
Our proof uses the methods of reducibility and discharging.  First, we prove
that certain properties must hold for every minimal counterexample to the
theorem (by ``minimal'' we mean having the fewest vertices and,
subject to that, the fewest non-triangular faces).  To conclude,
we give a counting argument, via the discharging method, showing that
every planar graph violates one of these properties.  Thus, no minimal
counterexample exists, so the theorem is true.  

Hereafter, we write $G$ to denote a minimal counterexample to the theorem.  To
remind the reader of this assumption, we will often refer to a \emph{minimal
$G$}.  Whenever we say ``a coloring'', we mean a $2$-fold $9$-coloring.
Note that $G$ is a plane triangulation; otherwise, adding an edge contradicts our
choice of $G$ as having the fewest non-triangular faces.

\begin{lem}\label{noTriangles}
A minimal $G$ has no separating clique.  Specifically, $G$ has no separating $3$-cycle.
\end{lem}
\begin{proof}
Suppose $G$ has a separating clique $X$ and let $C_1, \ldots, C_k$ be the components of
$G\setminus X$.  By minimality of $|G|$, we have colorings of $G[V(C_i) \cup X]$ for
each $i\in\set{1,\ldots,k}$.  Permute the colors on each subgraph
$G[V(C_i) \cup X]$ so the colorings agree on $X$.  Now identifying the copies of $X$
in each $G[V(C_i) \cup X]$ gives a coloring of $G$, a contradiction.
\end{proof}

Although it was easy to prove, Lemma~\ref{noTriangles} will play a crucial role
in our proof.  We will often want to identify two neighbors $u_1$ and $u_2$ of
a vertex $v$ and color the smaller graph by minimality.  To do so, we must
ensure that $u_1$ and $u_2$ are indeed non-adjacent; these arguments typically
use the fact that if $u_1$ and $u_2$ were adjacent, then $u_1u_2v$ would be a
separating 3-cycle.

\begin{lem}
A minimal $G$ has minimum degree $5$.
\end{lem}
\begin{proof}
Since $G$ is a plane triangulation, it has minimum degree at least $3$ and at most
$5$.  If $G$ contains a $3$-vertex, then its neighbors induce a separating $3$-cycle,
contradicting Lemma~\ref{noTriangles}.  If $G$ contains a $4$-vertex $v$,
then some pair of its neighbors are non-adjacent, since $K_5$ is non-planar. 
Form $G'$ from $G$ by deleting $v$ and contracting a non-adjacent pair of its
neighbors.  Color $G'$ by minimality, then lift the coloring back to $G$; only
$v$ is uncolored.  Since two of $v$'s neighbors have the same colors, we can extend
the coloring to $G$.
\end{proof}

The following fact will often allow us to extend a 2-fold 9-coloring to the uncolored
vertices of an induced $K_{1,3}$.  It will be useful in verifying that numerous
configurations are forbidden from a minimal $G$.
We will also often apply it when the uncolored subgraph is simply $P_3$.

\begin{fact}
Let $H=K_{1,3}$.  If each leaf has a list of size $3$ and the center vertex has a
list of size $5$, then we can choose $2$ colors for each vertex from its lists such
that adjacent vertices get disjoint sets of colors.
\label{fact1}
\end{fact}
\begin{proof}
Let $v$ denote the center vertex and $u_1,u_2,u_3$ the leaves.  
Since $2|L(v)| > |L(u_1)|+|L(u_2)|+|L(u_3)|$, some color $c\in L(v)$ appears in
$L(u_i)$ for at most one $u_i$. If such a $u_i$ exists, then by symmetry, say
it is $u_1$; now color $v$ with $c$ and some color not in $L(u_1)$.  Otherwise
color $v$ with $c$ and an arbitrary color.  Now color each $u_i$ arbitrarily
from its at least $2$ available colors.
\end{proof}

We use the same approach to prove each of Lemmas~\ref{556}, \ref{666}, and \ref{7566}.  
Our idea is to contract some edges of $G$ to get a smaller planar graph $G'$,
which we color by minimality.  In particular, in $G'$ we identify some pairs of
non-adjacent vertices of $G$ that each have a common neighbor.  When we lift the
coloring of $G'$ to $G$ this means that some of the uncolored vertices will
have neighbors with both colors the same, reducing the number of colors used on
the neighborhood of each such uncolored vertex.  

One early example of this technique is Kainen's proof~\cite{Kainen} of the
$5$ Color Theorem.
If $G$ is a planar graph, then by Euler's Theorem, $G$ has a $5^-$-vertex $v$. 
If $d(v)\le 4$, then we 5-color $G-v$ by minimality; now, since $d(v)\le 4$, we can
extend the 5-coloring to $v$.  Suppose instead that $d(v)=5$.  Since $K_6$ is
non-planar, $v$ has two neighbors $u_1$ and $u_2$ that are non-adjacent; form
$G'$ by contracting the edges $vu_1$ and $vu_2$, and again 5-color $G'$ by minimality.
To extend the 5-coloring to $v$, we note that even though $d(v)=5$, at most four
colors appear on the neighbors of $v$ (since $u_1$ and $u_2$ have the same
color).  This completes the proof.

Because a minimal $G$ has no separating 3-cycles, if vertices $u_1$ and $u_2$ have a
common neighbor $v$ and do not appear sequentially on the cycle induced by the
neighborhood of $v$, then $u_1$ and $u_2$ are non-adjacent.
The numeric labels in the figures denote pairs (or more) of vertices that are
identified in $G'$ when we delete any vertices labeled $v$, $u_1$, $u_2$ or $u_3$; vertices
with the same numeric label get identified.

Typically, it suffices to verify that the vertices receiving a common numeric label are
pairwise nonadjacent.  One potential complication is if two vertices that are drawn
as distinct are in fact the same vertex.  This usually cannot happen if the 
vertices have a common neighbor $v$, since then the degree of $v$ would be too small.
Similarly, it cannot happen if they are joined by a path of length three, 
since then we would get a separating 3-cycle.  

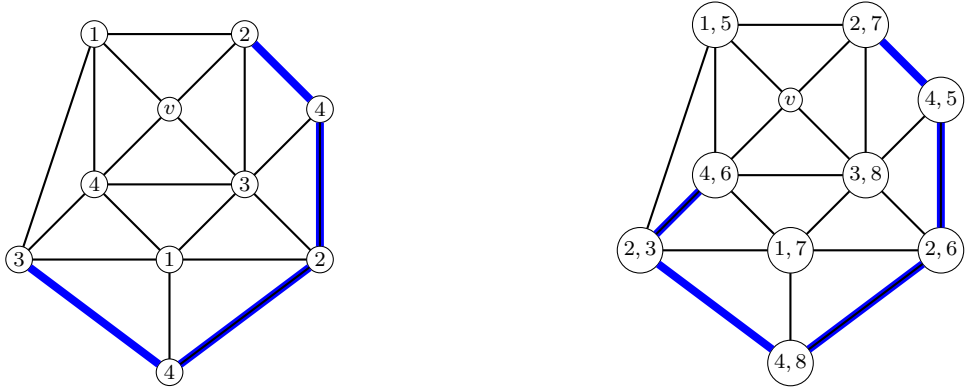
\begin{figure}[h!]
\subfloat[The $2,4$-path is blocked by the $1,3$-path.]{\makebox[.5\textwidth]{
\begin{tikzpicture}[scale = 10]
\tikzstyle{VertexStyle} = []
\tikzstyle{EdgeStyle} = []
\tikzstyle{labeledStyle}=[shape = circle, minimum size = 6pt, inner sep = 1.2pt, draw]
\tikzstyle{unlabeledStyle}=[shape = circle, minimum size = 6pt, inner sep = 1.2pt, draw, fill]
\tikzstyle{38cfed53-dd15-4981-97c1-4749de192b1c}=[blue, line width=3]
\Vertex[style = labeledStyle, x = 0.699999988079071, y = 0.799999997019768, L = \tiny {$1$}]{v0}
\Vertex[style = labeledStyle, x = 0.899999976158142, y = 0.799999997019768, L = \tiny {$2$}]{v1}
\Vertex[style = labeledStyle, x = 0.699999988079071, y = 0.599999994039536, L = \tiny {$4$}]{v2}
\Vertex[style = labeledStyle, x = 0.899999976158142, y = 0.599999994039536, L = \tiny {$3$}]{v3}
%\Vertex[style = labeledStyle, x = 0.600000023841858, y = 0.699999988079071, L = \tiny {$3$}]{v4}
\Vertex[style = labeledStyle, x = 0.600000023841858, y = 0.5, L = \tiny {$3$}]{v5}
\Vertex[style = labeledStyle, x = 0.800000011920929, y = 0.5, L = \tiny {$1$}]{v6}
\Vertex[style = labeledStyle, x = 1, y = 0.5, L = \tiny {$2$}]{v7}
\Vertex[style = labeledStyle, x = 1, y = 0.699999988079071, L = \tiny {$4$}]{v8}
%\Vertex[style = labeledStyle, x = 1.10000002384186, y = 0.599999994039536, L = \tiny {$2$}]{v9}
\Vertex[style = labeledStyle, x = .8, y = 0.350000023841858, L = \tiny {$4$}]{v10}
\Vertex[style = labeledStyle, x = 0.800000011920929, y = 0.699999988079071, L = \tiny {$v$}]{v11}
\Edge[](v0)(v5)
\Edge[](v1)(v0)
\Edge[](v2)(v0)
\Edge[](v1)(v3)
\Edge[](v2)(v3)
%\Edge[](v4)(v0)
%\Edge[](v4)(v5)
\Edge[](v5)(v6)
\Edge[](v7)(v3)
\Edge[](v7)(v6)
\Edge[style = 38cfed53-dd15-4981-97c1-4749de192b1c](v8)(v1)
\Edge[style = 38cfed53-dd15-4981-97c1-4749de192b1c](v8)(v7)
\Edge[style = 38cfed53-dd15-4981-97c1-4749de192b1c](v10)(v7)
\Edge[](v0)(v11)
\Edge[](v1)(v11)
\Edge[](v2)(v11)
\Edge[](v3)(v11)
\Edge[](v2)(v5)
\Edge[](v2)(v6)
\Edge[](v3)(v6)
%\Edge[](v4)(v2)
\Edge[](v8)(v3)
\Edge[](v8)(v7)
%\Edge[](v9)(v7)
\Edge[style = 38cfed53-dd15-4981-97c1-4749de192b1c](v10)(v5)
\Edge[](v10)(v6)
\Edge[](v10)(v7)
\end{tikzpicture}}}
\subfloat[The $2,4$-path gets through.]{\makebox[.5\textwidth]{
\begin{tikzpicture}[scale = 10]
\tikzstyle{VertexStyle} = []
\tikzstyle{EdgeStyle} = []
\tikzstyle{labeledStyle}=[shape = circle, minimum size = 6pt, inner sep = 1.2pt, draw]
\tikzstyle{unlabeledStyle}=[shape = circle, minimum size = 6pt, inner sep = 1.2pt, draw, fill]
\tikzstyle{38cfed53-dd15-4981-97c1-4749de192b1c}=[blue, line width=3]
\Vertex[style = labeledStyle, x = 0.699999988079071, y = 0.799999997019768, L =
\tiny {$1,5$}]{v0}
\Vertex[style = labeledStyle, x = 0.899999976158142, y = 0.799999997019768, L =
\tiny {$2,7$}]{v1}
\Vertex[style = labeledStyle, x = 0.699999988079071, y = 0.599999994039536, L =
\tiny {$4,6$}]{v2}
\Vertex[style = labeledStyle, x = 0.899999976158142, y = 0.599999994039536, L =
\tiny {$3,8$}]{v3}
%\Vertex[style = labeledStyle, x = 0.600000023841858, y = 0.699999988079071, L = \tiny {$3$}]{v4}
\Vertex[style = labeledStyle, x = 0.600000023841858, y = 0.5, L = \tiny {$2,3$}]{v5}
\Vertex[style = labeledStyle, x = 0.800000011920929, y = 0.5, L = \tiny {$1,7$}]{v6}
\Vertex[style = labeledStyle, x = 1, y = 0.5, L = \tiny {$2,6$}]{v7}
\Vertex[style = labeledStyle, x = 1, y = 0.699999988079071, L = \tiny {$4,5$}]{v8}
%\Vertex[style = labeledStyle, x = 1.10000002384186, y = 0.599999994039536, L = \tiny {$2$}]{v9}
\Vertex[style = labeledStyle, x = .8, y = 0.350000023841858, L = \tiny {$4,8$}]{v10}
\Vertex[style = labeledStyle, x = 0.800000011920929, y = 0.699999988079071, L = \tiny {$v$}]{v11}
\Edge[](v0)(v5)
\Edge[](v1)(v0)
\Edge[](v2)(v0)
\Edge[](v1)(v3)
\Edge[](v2)(v3)
%\Edge[](v4)(v0)
%\Edge[](v4)(v5)
\Edge[](v5)(v6)
\Edge[](v7)(v3)
\Edge[](v7)(v6)
\Edge[style = 38cfed53-dd15-4981-97c1-4749de192b1c](v2)(v5)
\Edge[style = 38cfed53-dd15-4981-97c1-4749de192b1c](v8)(v1)
\Edge[style = 38cfed53-dd15-4981-97c1-4749de192b1c](v8)(v7)
\Edge[style = 38cfed53-dd15-4981-97c1-4749de192b1c](v10)(v7)
\Edge[](v0)(v11)
\Edge[](v1)(v11)
\Edge[](v2)(v11)
\Edge[](v3)(v11)
\Edge[](v2)(v5)
\Edge[](v2)(v6)
\Edge[](v3)(v6)
%\Edge[](v4)(v2)
\Edge[](v8)(v3)
\Edge[](v8)(v7)
%\Edge[](v9)(v7)
\Edge[style = 38cfed53-dd15-4981-97c1-4749de192b1c](v10)(v5)
\Edge[](v10)(v6)
\Edge[](v10)(v7)
\end{tikzpicture}}}
\caption{The problem with Kempe chains for $2$-fold coloring.}
\label{badkempe}
\end{figure}

For $4$-coloring, Birkhoff \cite{birkhoff} showed how to exclude separating $4$-cycles and $5$-cycles.  
Excluding separating $4$-cycles would simplify our arguments below since we
would not need to worry about vertices at distance at most four being the same.
The proof excluding $4$-cycles for $4$-coloring is quite easy, but it does not
work in our context because standard Kempe chain arguments break down for $2$-fold coloring.
The problem is illustrated in Figure~\ref{badkempe}. Figure~\ref{badkempe}(A)
shows the situation for $1$-fold coloring; here the $13$-path blocks the
$24$-path. Figure~\ref{badkempe}(B) shows the situation for $2$-fold coloring,
here the $24$-path can get through because on the $13$-path, a vertex has color
$2$ as well as color $1$.
%\documentclass{article}
%\usepackage{tikz,tkz-graph}
%\usepackage{subfig}
%\begin{document}
%\end{document}

\begin{lem}
A minimal $G$ has no $5$-vertex with a $5$-neighbor and a non-adjacent $6^-$-neighbor.
\label{556}
\end{lem}
\begin{proof}
We first consider the case where a $5$-vertex $v$ has non-adjacent $5$-neighbors $u_1$
and $u_2$, as shown in Figure~\ref{fig556}(A). We color $G'$ by minimality, then lift the
coloring to $G$.  (Recall that to form $G'$, we delete $v$ and all $u_i$ and for
each pair (or more) of vertices with the same label, we identify them.)
Now in $G$ each $u_i$ has a list of at least $3$ colors and
$v$ has a list of at least $5$ colors; so, by Fact~\ref{fact1}, we can extend the coloring to $G$.
%\url{http://tinyurl.com/msl9eok}.  
\begin{figure}[ht!]
%\centering
\subfloat[A 5-vertex, $v$, with non-adjacent 5-neighbors, $u_1$ and $u_2$.]{\makebox[.5\textwidth]{
\begin{tikzpicture}[scale = 9]
\tikzstyle{VertexStyle}=[shape = circle, minimum size = 6pt, inner sep = 1.2pt, draw]
\Vertex[x = 0.60, y = 0.75, L = \tiny {$v$}]{v0}
\Vertex[x = 0.45, y = 0.80, L = \tiny {$u_1$}]{v1}
\Vertex[x = 0.60, y = 0.95, L = \tiny {$2$}]{v2}
\Vertex[x = 0.50, y = 0.60, L = \tiny {$1$}]{v3}
\Vertex[x = 0.70, y = 0.60, L = \tiny {$2$}]{v4}
\Vertex[x = 0.75, y = 0.80, L = \tiny {$u_2$}]{v5}
\Vertex[x = 0.30, y = 0.85, L = \tiny {$1$}]{v6}
\Vertex[x = 0.30, y = 0.65, L = \tiny {}]{v7}
\Vertex[x = 0.90, y = 0.85, L = \tiny {}]{v8}
\Vertex[x = 0.90, y = 0.65, L = \tiny {}]{v9}
\Edge[](v0)(v1)
\Edge[](v0)(v3)
\Edge[](v0)(v4)
\Edge[](v0)(v5)
\Edge[](v0)(v2)
\Edge[](v2)(v8)
\Edge[](v2)(v5)
\Edge[](v5)(v4)
\Edge[](v5)(v9)
\Edge[](v5)(v8)
\Edge[](v9)(v8)
\Edge[](v9)(v4)
\Edge[](v1)(v6)
\Edge[](v1)(v7)
\Edge[](v1)(v3)
\Edge[](v3)(v4)
\Edge[](v3)(v7)
\Edge[](v6)(v7)
\Edge[](v2)(v1)
\Edge[](v2)(v6)
\end{tikzpicture}}}
\subfloat[A 5-vertex, $v$, with a non-adjacent 5-neighbor, $u_1$, and
6-neighbor, $u_2$.]{\makebox[.5\textwidth]{
\begin{tikzpicture}[scale = 9]
\tikzstyle{VertexStyle}=[shape = circle, minimum size = 6pt, inner sep = 1.2pt, draw]
\Vertex[x = 0.60, y = 0.75, L = \tiny {$v$}]{v0}
\Vertex[x = 0.45, y = 0.80, L = \tiny {$u_1$}]{v1}
\Vertex[x = 0.60, y = 0.95, L = \tiny {$2$}]{v2}
\Vertex[x = 0.50, y = 0.60, L = \tiny {$1$}]{v3}
\Vertex[x = 0.70, y = 0.60, L = \tiny {$2$}]{v4}
\Vertex[x = 0.80, y = 0.75, L = \tiny {$u_2$}]{v5}
\Vertex[x = 0.30, y = 0.85, L = \tiny {$1$}]{v6}
\Vertex[x = 0.30, y = 0.65, L = \tiny {}]{v7}
\Vertex[x = 1.00, y = 0.75, L = \tiny {}]{v8}
\Vertex[x = 0.90, y = 0.60, L = \tiny {$3$}]{v9}
\Vertex[x = 0.90, y = 0.90, L = \tiny {$3$}]{v10}
\Edge[](v0)(v1)
\Edge[](v0)(v3)
\Edge[](v0)(v4)
\Edge[](v0)(v5)
\Edge[](v0)(v2)
\Edge[](v2)(v5)
\Edge[](v5)(v4)
\Edge[](v5)(v9)
\Edge[](v5)(v8)
\Edge[](v9)(v8)
\Edge[](v9)(v4)
\Edge[](v1)(v6)
\Edge[](v1)(v7)
\Edge[](v1)(v3)
\Edge[](v3)(v4)
\Edge[](v3)(v7)
\Edge[](v6)(v7)
\Edge[](v2)(v1)
\Edge[](v2)(v6)
\Edge[](v10)(v8)
\Edge[](v10)(v5)
\Edge[](v10)(v2)
\end{tikzpicture}}}
\caption{The cases of Lemma~\ref{556}.}
\label{fig556}
\end{figure}
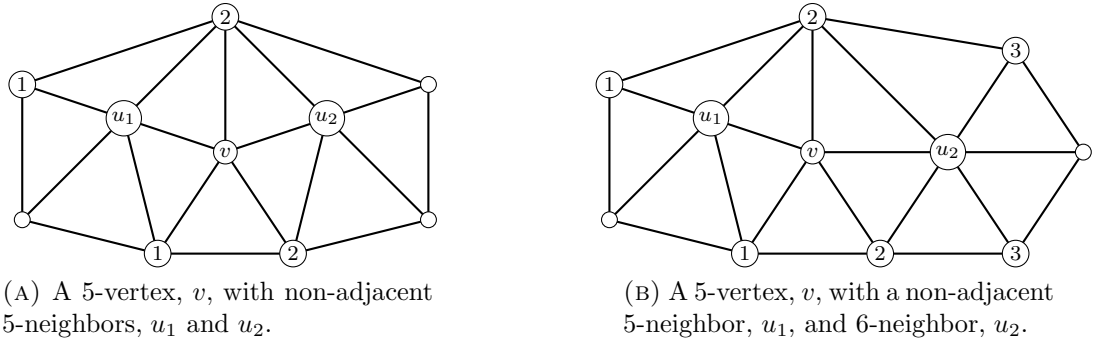

Now we consider the case where a $5$-vertex $v$ has a $5$-neighbor and a $6$-neighbor
that are non-adjacent, as shown in Figure~\ref{fig556}(B). 
%\url{http://tinyurl.com/lh48o9n}.  
Again, when we lift the coloring of $G'$ to $G$, $v$ has a list of size $5$ and
each of its uncolored neighbors
has a list of size $3$.  Hence, by Fact~\ref{fact1}, we can extend the
coloring of $G'$ to $G$.  
Here no pair of labeled vertices can be identified, since each such pair is
drawn at distance three or less (and $G$ has no separating 3-cycle).
\end{proof}

\begin{lem}
A minimal $G$ has no $6$-vertex with non-adjacent $6^-$-neighbors.
\label{666}
\end{lem}
\begin{proof}
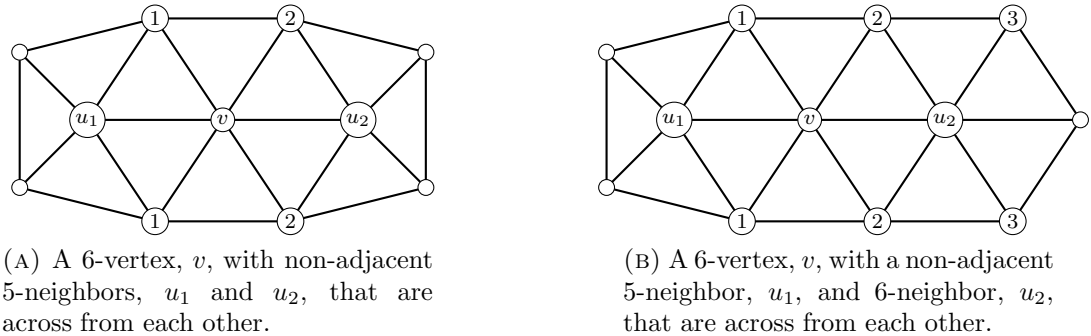
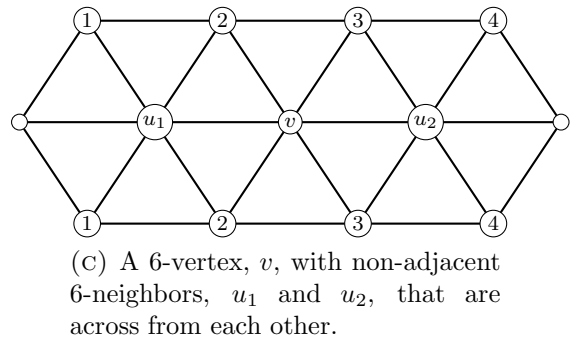
\begin{figure}[!ht]
%\centering
\subfloat[A 6-vertex, $v$, with non-adjacent 5-neighbors, $u_1$ and $u_2$, that
are across from each other.]{\makebox[.5\textwidth]{
\begin{tikzpicture}[scale = 9]
\tikzstyle{VertexStyle}=[shape = circle, minimum size = 6pt, inner sep = 1.2pt, draw]
\Vertex[x = 0.60, y = 0.75, L = \tiny {$v$}]{v0}
\Vertex[x = 0.50, y = 0.90, L = \tiny {$1$}]{v1}
\Vertex[x = 0.40, y = 0.75, L = \tiny {$u_1$}]{v2}
\Vertex[x = 0.70, y = 0.90, L = \tiny {$2$}]{v3}
\Vertex[x = 0.50, y = 0.60, L = \tiny {$1$}]{v4}
\Vertex[x = 0.70, y = 0.60, L = \tiny {$2$}]{v5}
\Vertex[x = 0.80, y = 0.75, L = \tiny {$u_2$}]{v6}
\Vertex[x = 0.30, y = 0.85, L = \tiny {}]{v7}
\Vertex[x = 0.30, y = 0.65, L = \tiny {}]{v8}
\Vertex[x = 0.90, y = 0.85, L = \tiny {}]{v9}
\Vertex[x = 0.90, y = 0.65, L = \tiny {}]{v10}
\Edge[](v0)(v1)
\Edge[](v0)(v2)
\Edge[](v0)(v4)
\Edge[](v0)(v5)
\Edge[](v0)(v6)
\Edge[](v0)(v3)
\Edge[](v3)(v1)
\Edge[](v3)(v9)
\Edge[](v3)(v6)
\Edge[](v6)(v5)
\Edge[](v6)(v10)
\Edge[](v6)(v9)
\Edge[](v10)(v9)
\Edge[](v10)(v5)
\Edge[](v2)(v1)
\Edge[](v2)(v7)
\Edge[](v2)(v8)
\Edge[](v2)(v4)
\Edge[](v4)(v5)
\Edge[](v4)(v8)
\Edge[](v7)(v1)
\Edge[](v7)(v8)
\end{tikzpicture}}}
%%%%
\subfloat[A 6-vertex, $v$, with a non-adjacent 5-neighbor, $u_1$, and
6-neighbor, $u_2$, that are across from each other.]{\makebox[.5\textwidth]{
\begin{tikzpicture}[scale = 9]
\tikzstyle{VertexStyle}=[shape = circle, minimum size = 6pt, inner sep = 1.2pt, draw]
\Vertex[x = 0.60, y = 0.75, L = \tiny {$v$}]{v0}
\Vertex[x = 0.50, y = 0.90, L = \tiny {$1$}]{v1}
\Vertex[x = 0.40, y = 0.75, L = \tiny {$u_1$}]{v2}
\Vertex[x = 0.70, y = 0.90, L = \tiny {$2$}]{v3}
\Vertex[x = 0.50, y = 0.60, L = \tiny {$1$}]{v4}
\Vertex[x = 0.70, y = 0.60, L = \tiny {$2$}]{v5}
\Vertex[x = 0.80, y = 0.75, L = \tiny {$u_2$}]{v6}
\Vertex[x = 0.30, y = 0.85, L = \tiny {}]{v7}
\Vertex[x = 0.30, y = 0.65, L = \tiny {}]{v8}
\Vertex[x = 0.90, y = 0.90, L = \tiny {$3$}]{v9}
\Vertex[x = 0.90, y = 0.60, L = \tiny {$3$}]{v10}
\Vertex[x = 1, y = 0.75, L = \tiny {}]{v11}
\Edge[](v0)(v1)
\Edge[](v0)(v2)
\Edge[](v0)(v4)
\Edge[](v0)(v5)
\Edge[](v0)(v6)
\Edge[](v0)(v3)
\Edge[](v3)(v1)
\Edge[](v3)(v9)
\Edge[](v3)(v6)
\Edge[](v6)(v5)
\Edge[](v6)(v10)
\Edge[](v6)(v9)
\Edge[](v10)(v5)
\Edge[](v2)(v1)
\Edge[](v2)(v7)
\Edge[](v2)(v8)
\Edge[](v2)(v4)
\Edge[](v4)(v5)
\Edge[](v4)(v8)
\Edge[](v7)(v1)
\Edge[](v7)(v8)
\Edge[](v9)(v11)
\Edge[](v11)(v6)
\Edge[](v11)(v10)
\end{tikzpicture}}}
  %\caption{A subfigure}
%\end{figure}
%
%\subfloat[A 5-vertex, $v$, with a non-adjacent 5-neighbor, $u_1$, and
%6-neighbor, $u_2$.]{\makebox[.5\textwidth]{

\subfloat[A 6-vertex, $v$, with non-adjacent 6-neighbors, $u_1$ and $u_2$, that
are across from each other.]{\makebox[.5\textwidth]{
\begin{tikzpicture}[scale = 9]
\tikzstyle{VertexStyle}=[shape = circle, minimum size = 6pt, inner sep = 1.2pt, draw]
\Vertex[x = 0.60, y = 0.75, L = \tiny {$v$}]{v0}
\Vertex[x = 0.50, y = 0.90, L = \tiny {$2$}]{v1}
\Vertex[x = 0.40, y = 0.75, L = \tiny {$u_1$}]{v2}
\Vertex[x = 0.70, y = 0.90, L = \tiny {$3$}]{v3}
\Vertex[x = 0.50, y = 0.60, L = \tiny {$2$}]{v4}
\Vertex[x = 0.70, y = 0.60, L = \tiny {$3$}]{v5}
\Vertex[x = 0.80, y = 0.75, L = \tiny {$u_2$}]{v6}
\Vertex[x = 0.30, y = 0.90, L = \tiny {$1$}]{v7}
\Vertex[x = 0.30, y = 0.60, L = \tiny {$1$}]{v8}
\Vertex[x = 0.90, y = 0.90, L = \tiny {$4$}]{v9}
\Vertex[x = 0.90, y = 0.60, L = \tiny {$4$}]{v10}
\Vertex[x = 1.00, y = 0.75, L = \tiny {}]{v11}
\Vertex[x = 0.20, y = 0.75, L = \tiny {}]{v12}
\Edge[](v0)(v1)
\Edge[](v0)(v2)
\Edge[](v0)(v4)
\Edge[](v0)(v5)
\Edge[](v0)(v6)
\Edge[](v0)(v3)
\Edge[](v3)(v1)
\Edge[](v3)(v9)
\Edge[](v3)(v6)
\Edge[](v6)(v5)
\Edge[](v6)(v10)
\Edge[](v6)(v9)
\Edge[](v10)(v5)
\Edge[](v2)(v1)
\Edge[](v2)(v7)
\Edge[](v2)(v8)
\Edge[](v2)(v4)
\Edge[](v4)(v5)
\Edge[](v4)(v8)
\Edge[](v7)(v1)
\Edge[](v9)(v11)
\Edge[](v11)(v6)
\Edge[](v11)(v10)
\Edge[](v12)(v2)
\Edge[](v12)(v7)
\Edge[](v12)(v8)
\end{tikzpicture}}}
\caption{The ``across'' cases of Lemma~\ref{666}.}
\label{fig666across}
\end{figure}
Let $v$ be a $6$-vertex with two non-adjacent $6^-$-neighbors, $u_1$ and $u_2$.
We have three possibilities for the degrees of these $6^-$-neighbors: two
$5$-vertices, a $5$-vertex and a $6$-vertex, and two $6$-vertices.
For each choice of degrees for the $u_i$s, we have two possibilities for their relative
location; they could be ``across'' from each other (at distance three along the
cycle induced by the neighbors of $v$) or ``offset'' from each other (at
distance two along the same cycle).  This yields a total of six possibilities; the
three across possibilities are shown in Figure~\ref{fig666across} and the
three offset possibilities are shown in Figure~\ref{fig666offset}.

In Figures~\ref{fig666across}(A,B), all of the vertices with numeric
labels (those that will be identified in $G'$) must be distinct, since they are
drawn within distance three of each other.
%\url{http://tinyurl.com/n545sqf}; 566 (across): http://tinyurl.com/on2skct;
%666 (across): http://tinyurl.com/oy94ja7. 
%
The only complication is in Figure~\ref{fig666across}(C): a
vertex labeled $1$ might be the same as a vertex labeled $4$ that is drawn at distance
four; call this vertex $x$.  By symmetry, assume that $x$ is formed by
identifying the vertex in the top left labeled 1 and the vertex in the bottom
right labeled 4.  This is only a problem if also a vertex
labeled $1$ is adjacent to one labeled $4$; so suppose this happens.  
Note that the vertex in the top right labeled 4 cannot be adjacent to the vertex
in the bottom left labeled 1; they are on opposite sides of the cycle
$xu_1vu_2$.  So, again by symmetry, we assume that $x$ is adjacent to the vertex
in the bottom left labeled 1.  However, now we have a separating 3-cycle
(consisting of $x$, its neighbor labeled 1, and their common neighbor $u_1$);
this contradicts Lemma~\ref{noTriangles}.  This contradiction finishes the across cases.
%
%Now here are pictures for the ``offset'' case: 
%565 (offset): \url{http://tinyurl.com/lomkluw};
%566 (offset): \url{http://tinyurl.com/p82ryrg};
%666 (offset): \url{http://tinyurl.com/ohuehaq}.
\begin{figure}[!ht]
%\centering
\subfloat[A 6-vertex, $v$, with non-adjacent 5-neighbors, $u_1$ and $u_2$, that
are offset from each other.]{\makebox[.5\textwidth]{
\begin{tikzpicture}[scale = 9]
\tikzstyle{VertexStyle}=[shape = circle, minimum size = 6pt, inner sep = 1.2pt, draw]
\Vertex[x = 0.80, y = 0.60, L = \tiny {$v$}]{v0}
\Vertex[x = 0.70, y = 0.75, L = \tiny {$u_1$}]{v1}
\Vertex[x = 0.60, y = 0.60, L = \tiny {$1$}]{v2}
\Vertex[x = 0.90, y = 0.75, L = \tiny {$1$}]{v3}
\Vertex[x = 0.70, y = 0.45, L = \tiny {}]{v4}
\Vertex[x = 0.90, y = 0.45, L = \tiny {$1$}]{v5}
\Vertex[x = 1.00, y = 0.60, L = \tiny {$u_2$}]{v6}
\Vertex[x = 0.50, y = 0.75, L = \tiny {}]{v7}
\Vertex[x = 1.10, y = 0.75, L = \tiny {}]{v8}
\Vertex[x = 1.10, y = 0.45, L = \tiny {}]{v9}
\Vertex[x = 0.80, y = 0.90, L = \tiny {}]{v10}
\Edge[](v0)(v1)
\Edge[](v0)(v2)
\Edge[](v0)(v4)
\Edge[](v0)(v5)
\Edge[](v0)(v6)
\Edge[](v0)(v3)
\Edge[](v3)(v1)
\Edge[](v3)(v8)
\Edge[](v3)(v6)
\Edge[](v6)(v5)
\Edge[](v6)(v9)
\Edge[](v6)(v8)
\Edge[](v9)(v5)
\Edge[](v2)(v1)
\Edge[](v2)(v7)
\Edge[](v2)(v4)
\Edge[](v4)(v5)
\Edge[](v7)(v1)
\Edge[](v10)(v1)
\Edge[](v10)(v3)
\Edge[](v8)(v9)
\Edge[](v10)(v7)
\end{tikzpicture}}}
%%%%
\subfloat[A 6-vertex, $v$, with a non-adjacent 5-neighbor, $u_1$, and
6-neighbor, $u_2$, that are offset from each other.]{\makebox[.5\textwidth]{
\begin{tikzpicture}[scale = 9]
\tikzstyle{VertexStyle}=[shape = circle, minimum size = 6pt, inner sep = 1.2pt, draw]
\Vertex[x = 0.80, y = 0.60, L = \tiny {$v$}]{v0}
\Vertex[x = 0.70, y = 0.75, L = \tiny {$u_1$}]{v1}
\Vertex[x = 0.60, y = 0.60, L = \tiny {$1$}]{v2}
\Vertex[x = 0.90, y = 0.75, L = \tiny {$1$}]{v3}
\Vertex[x = 0.70, y = 0.45, L = \tiny {}]{v4}
\Vertex[x = 0.90, y = 0.45, L = \tiny {$1$}]{v5}
\Vertex[x = 1.00, y = 0.60, L = \tiny {$u_2$}]{v6}
\Vertex[x = 0.50, y = 0.75, L = \tiny {}]{v7}
\Vertex[x = 1.10, y = 0.75, L = \tiny {$2$}]{v8}
\Vertex[x = 1.10, y = 0.45, L = \tiny {$2$}]{v9}
\Vertex[x = 0.80, y = 0.90, L = \tiny {}]{v10}
\Vertex[x = 1.20, y = 0.60, L = \tiny {}]{v11}
\Edge[](v0)(v1)
\Edge[](v0)(v2)
\Edge[](v0)(v4)
\Edge[](v0)(v5)
\Edge[](v0)(v6)
\Edge[](v0)(v3)
\Edge[](v3)(v1)
\Edge[](v3)(v8)
\Edge[](v3)(v6)
\Edge[](v6)(v5)
\Edge[](v6)(v9)
\Edge[](v6)(v8)
\Edge[](v9)(v5)
\Edge[](v2)(v1)
\Edge[](v2)(v7)
\Edge[](v2)(v4)
\Edge[](v4)(v5)
\Edge[](v7)(v1)
\Edge[](v10)(v1)
\Edge[](v10)(v3)
\Edge[](v10)(v7)
\Edge[](v11)(v6)
\Edge[](v11)(v8)
\Edge[](v11)(v9)
\end{tikzpicture}}}

\subfloat[A 6-vertex, $v$, with non-adjacent 6-neighbors, $u_1$ and
$u_2$, that are offset from each other (case i).]{\makebox[.5\textwidth]{
\begin{tikzpicture}[scale = 9]
\tikzstyle{VertexStyle}=[shape = circle, minimum size = 6pt, inner sep = 1.2pt, draw]
\Vertex[x = 0.80, y = 0.60, L = \tiny {$v$}]{v0}
\Vertex[x = 0.70, y = 0.75, L = \tiny {$u_1$}]{v1}
\Vertex[x = 0.60, y = 0.60, L = \tiny {$2$}]{v2}
\Vertex[x = 0.90, y = 0.75, L = \tiny {$2$}]{v3}
\Vertex[x = 0.70, y = 0.45, L = \tiny {}]{v4}
\Vertex[x = 0.90, y = 0.45, L = \tiny {$2$}]{v5}
\Vertex[x = 1.00, y = 0.60, L = \tiny {$u_2$}]{v6}
\Vertex[x = 0.50, y = 0.75, L = \tiny {$1$}]{v7}
\Vertex[x = 1.10, y = 0.75, L = \tiny {$3$}]{v8}
\Vertex[x = 1.10, y = 0.45, L = \tiny {$3$}]{v9}
\Vertex[x = 0.80, y = 0.90, L = \tiny {$1$}]{v10}
\Vertex[x = 1.20, y = 0.60, L = \tiny {}]{v11}
\Vertex[x = 0.60, y = 0.90, L = \tiny {}]{v12}
\Edge[](v0)(v1)
\Edge[](v0)(v2)
\Edge[](v0)(v4)
\Edge[](v0)(v5)
\Edge[](v0)(v6)
\Edge[](v0)(v3)
\Edge[](v3)(v1)
\Edge[](v3)(v8)
\Edge[](v3)(v6)
\Edge[](v6)(v5)
\Edge[](v6)(v9)
\Edge[](v6)(v8)
\Edge[](v9)(v5)
\Edge[](v2)(v1)
\Edge[](v2)(v7)
\Edge[](v2)(v4)
\Edge[](v4)(v5)
\Edge[](v7)(v1)
\Edge[](v10)(v1)
\Edge[](v10)(v3)
\Edge[](v11)(v6)
\Edge[](v11)(v8)
\Edge[](v11)(v9)
\Edge[](v12)(v1)
\Edge[](v12)(v7)
\Edge[](v12)(v10)
\end{tikzpicture}}}
%%%%
\begin{comment}
\subfloat[A 5-vertex, $v$, with a non-adjacent 5-neighbor, $u_1$, and
6-neighbor, $u_2$.]{\makebox[.5\textwidth]{
\begin{tikzpicture}[scale = 9]
\tikzstyle{VertexStyle}=[shape = circle, minimum size = 6pt, inner sep = 1.2pt, draw]
\Vertex[x = 0.55, y = 0.55, L = \tiny {$v$}]{v0}
\Vertex[x = 0.45, y = 0.70, L = \tiny {$u_1$}]{v1}
\Vertex[x = 0.35, y = 0.55, L = \tiny {$1$}]{v2}
\Vertex[x = 0.65, y = 0.70, L = \tiny {$1$}]{v3}
\Vertex[x = 0.45, y = 0.40, L = \tiny {}]{v4}
\Vertex[x = 0.65, y = 0.40, L = \tiny {$1$}]{v5}
\Vertex[x = 0.75, y = 0.55, L = \tiny {$u_2$}]{v6}
\Vertex[x = 0.85, y = 0.70, L = \tiny {}]{v7}
\Vertex[x = 0.55, y = 0.85, L = \tiny {}]{v8}
\Vertex[x = 0.95, y = 0.55, L = \tiny {$1$}]{v9}
\Vertex[x = 0.35, y = 0.85, L = \tiny {$1$}]{v10}
\Vertex[x = 0.85, y = 0.40, L = \tiny {}]{v11}
\Vertex[x = 0.25, y = 0.70, L = \tiny {}]{v12}
\Edge[](v0)(v1)
\Edge[](v0)(v2)
\Edge[](v0)(v4)
\Edge[](v0)(v5)
\Edge[](v0)(v6)
\Edge[](v0)(v3)
\Edge[](v3)(v1)
\Edge[](v3)(v7)
\Edge[](v3)(v6)
\Edge[](v6)(v5)
\Edge[](v6)(v7)
\Edge[](v2)(v1)
\Edge[](v2)(v4)
\Edge[](v4)(v5)
\Edge[](v8)(v1)
\Edge[](v8)(v3)
\Edge[](v9)(v7)
\Edge[](v9)(v6)
\Edge[](v10)(v8)
\Edge[](v10)(v1)
\Edge[](v8)(v7)
\Edge[](v11)(v9)
\Edge[](v11)(v6)
\Edge[](v11)(v5)
\Edge[](v12)(v2)
\Edge[](v12)(v10)
\Edge[](v12)(v1)
\end{tikzpicture}}}
%
%
\end{comment}
\subfloat[A 6-vertex, $v$, with non-adjacent 6-neighbors, $u_1$ and
$u_2$, that are offset from each other (case ii).]{\makebox[.5\textwidth]{
\begin{tikzpicture}[scale = 9]
\tikzstyle{VertexStyle}=[shape = circle, minimum size = 6pt, inner sep = 1.2pt, draw]
\Vertex[x = 0.55, y = 0.55, L = \tiny {$v$}]{v0}
\Vertex[x = 0.45, y = 0.70, L = \tiny {$u_1$}]{v1}
\Vertex[x = 0.35, y = 0.55, L = \tiny {$1$}]{v2}
\Vertex[x = 0.65, y = 0.70, L = \tiny {$1$}]{v3}
\Vertex[x = 0.45, y = 0.40, L = \tiny {}]{v4}
\Vertex[x = 0.65, y = 0.40, L = \tiny {$1$}]{v5}
\Vertex[x = 0.75, y = 0.55, L = \tiny {$u_2$}]{v6}
\Vertex[x = 0.85, y = 0.70, L = \tiny {$2$}]{v7}
\Vertex[x = 0.55, y = 0.85, L = \tiny {}]{v8}
\Vertex[x = 0.95, y = 0.55, L = \tiny {}]{v9}
\Vertex[x = 0.35, y = 0.85, L = \tiny {$1$}]{v10}
\Vertex[x = 0.85, y = 0.40, L = \tiny {$2$}]{v11}
%\Vertex[x = 0.25, y = 0.70, L = \tiny {}]{v12}
\Edge[](v0)(v1)
\Edge[](v0)(v2)
\Edge[](v0)(v4)
\Edge[](v0)(v5)
\Edge[](v0)(v6)
\Edge[](v0)(v3)
\Edge[](v3)(v1)
\Edge[](v3)(v7)
\Edge[](v3)(v6)
\Edge[](v6)(v5)
\Edge[](v6)(v7)
\Edge[](v2)(v1)
\Edge[](v2)(v4)
\Edge[](v4)(v5)
\Edge[](v8)(v1)
\Edge[](v8)(v3)
\Edge[](v9)(v7)
\Edge[](v9)(v6)
\Edge[](v10)(v8)
\Edge[](v10)(v1)
\Edge[](v11)(v9)
\Edge[](v11)(v6)
\Edge[](v11)(v5)
%\Edge[](v12)(v2)
%\Edge[](v12)(v1)
%\Edge[](v12)(v10)
\end{tikzpicture}}}

\caption{The ``offset'' cases of Lemma~\ref{666}.}
\label{fig666offset}
\end{figure}
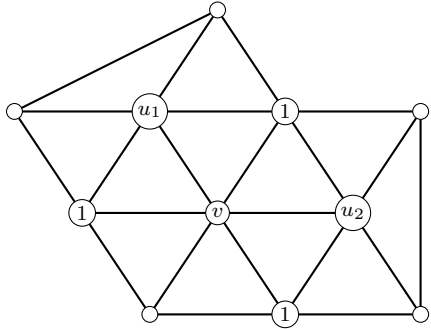
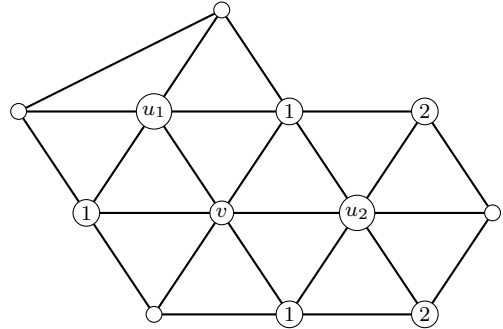
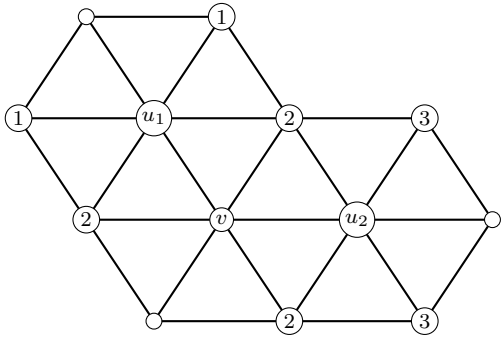
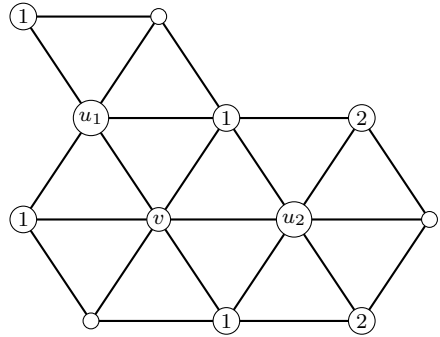

Now we consider the three offset cases, which are shown in Figure~\ref{fig666offset}.
As with the across cases, in Figures~\ref{fig666offset}(A,B), all of the
vertices with numeric labels must be distinct, since they are drawn within distance $3$
of each other.

The only complication in is the third case, shown in Figures~\ref{fig666offset}(C,D):
the vertices labeled $1$ and $3$ that are drawn at
distance four in Figure~\ref{fig666offset}(C) might be the same; if so, then
call this vertex $x$.  In this case we switch to the identifications shown in
Figure~\ref{fig666offset}(D); we omit from Figure~\ref{fig666offset}(D) a few
edges incident to $x$, to keep the picture pretty. Now all vertices with
numeric labels are at distance at most three, due to the extra edges incident
to $x$. Also, the two vertices labeled $1$ that are drawn at distance three are
non-adjacent, since they are separated by cycle $u_1vu_2x$.  This finishes the offset
cases.
\end{proof}

\begin{lem}
A minimal $G$ has no $7$-vertex with a $5$-neighbor and two other $6^-$-neighbors
such that all three are pairwise non-adjacent.
\label{7566}
\end{lem}
\begin{proof}
\begin{figure}[]
%\centering
\subfloat[A 7-vertex, $v$, with non-adjacent 5-neighbors, $u_1$, $u_2$, and
$u_3$.]{\makebox[.5\textwidth]{
\begin{tikzpicture}[scale = 10]
\tikzstyle{VertexStyle}=[shape = circle, minimum size = 6pt, inner sep = 1.2pt, draw]
\Vertex[x = 0.55, y = 0.55, L = \tiny {$v$}]{v0}
\Vertex[x = 0.40, y = 0.65, L = \tiny {$1$}]{v1}
\Vertex[x = 0.35, y = 0.50, L = \tiny {$u_1$}]{v2}
\Vertex[x = 0.70, y = 0.65, L = \tiny {$1$}]{v3}
\Vertex[x = 0.45, y = 0.40, L = \tiny {$1$}]{v4}
\Vertex[x = 0.65, y = 0.40, L = \tiny {$2$}]{v5}
\Vertex[x = 0.75, y = 0.50, L = \tiny {$u_3$}]{v6}
\Vertex[x = 0.55, y = 0.75, L = \tiny {$u_2$}]{v7}
\Vertex[x = 0.85, y = 0.40, L = \tiny {}]{v8}
\Vertex[x = 0.85, y = 0.60, L = \tiny {$2$}]{v9}
\Vertex[x = 0.25, y = 0.60, L = \tiny {}]{v10}
\Vertex[x = 0.25, y = 0.40, L = \tiny {}]{v11}
\Vertex[x = 0.40, y = 0.85, L = \tiny {}]{v13}
\Vertex[x = 0.70, y = 0.85, L = \tiny {}]{v14}
\Edge[](v0)(v1)
\Edge[](v0)(v2)
\Edge[](v0)(v4)
\Edge[](v0)(v5)
\Edge[](v0)(v6)
\Edge[](v0)(v3)
\Edge[](v3)(v6)
\Edge[](v6)(v5)
\Edge[](v2)(v1)
\Edge[](v2)(v4)
\Edge[](v4)(v5)
\Edge[](v7)(v0)
\Edge[](v7)(v3)
\Edge[](v7)(v1)
\Edge[](v9)(v6)
\Edge[](v9)(v8)
\Edge[](v8)(v5)
\Edge[](v8)(v6)
\Edge[](v9)(v3)
\Edge[](v11)(v10)
\Edge[](v11)(v2)
\Edge[](v11)(v4)
\Edge[](v10)(v2)
\Edge[](v10)(v1)
\Edge[](v13)(v14)
\Edge[](v13)(v7)
\Edge[](v14)(v7)
\Edge[](v1)(v13)
\Edge[](v3)(v14)
\end{tikzpicture}}}
~~~~~~
% Figure (B)
\subfloat[A 7-vertex, $v$, with a 6-neighbor, $u_3$, and
two 5-neighbors, $u_1$ and $u_2$, with all pairs of $u_i$s non-adjacent.]{\makebox[.5\textwidth]{
\begin{tikzpicture}[xscale = -10, yscale = 10]
\tikzstyle{VertexStyle}=[shape = circle, minimum size = 6pt, inner sep = 1.2pt, draw]
\Vertex[x = 0.55, y = 0.55, L = \tiny {$v$}]{v0}
\Vertex[x = 0.40, y = 0.65, L = \tiny {$2$}]{v1}
\Vertex[x = 0.35, y = 0.50, L = \tiny {$u_3$}]{v2}
\Vertex[x = 0.70, y = 0.65, L = \tiny {$2$}]{v3}
\Vertex[x = 0.45, y = 0.40, L = \tiny {$2$}]{v4}
\Vertex[x = 0.65, y = 0.40, L = \tiny {$1$}]{v5}
\Vertex[x = 0.75, y = 0.50, L = \tiny {$u_1$}]{v6}
\Vertex[x = 0.55, y = 0.75, L = \tiny {$u_2$}]{v7}
\Vertex[x = 0.85, y = 0.40, L = \tiny {}]{v8}
\Vertex[x = 0.85, y = 0.60, L = \tiny {$1$}]{v9}
\Vertex[x = 0.25, y = 0.60, L = \tiny {$3$}]{v10}
\Vertex[x = 0.25, y = 0.40, L = \tiny {$3$}]{v11}
\Vertex[x = 0.15, y = 0.50, L = \tiny {}]{v12}
\Vertex[x = 0.40, y = 0.85, L = \tiny {}]{v13}
\Vertex[x = 0.70, y = 0.85, L = \tiny {}]{v14}
\Edge[](v0)(v1)
\Edge[](v0)(v2)
\Edge[](v0)(v4)
\Edge[](v0)(v5)
\Edge[](v0)(v6)
\Edge[](v0)(v3)
\Edge[](v3)(v6)
\Edge[](v6)(v5)
\Edge[](v2)(v1)
\Edge[](v2)(v4)
\Edge[](v4)(v5)
\Edge[](v7)(v0)
\Edge[](v7)(v3)
\Edge[](v7)(v1)
\Edge[](v9)(v6)
\Edge[](v9)(v8)
\Edge[](v8)(v5)
\Edge[](v8)(v6)
\Edge[](v9)(v3)
\Edge[](v11)(v2)
\Edge[](v11)(v4)
\Edge[](v10)(v2)
\Edge[](v10)(v1)
\Edge[](v12)(v10)
\Edge[](v12)(v2)
\Edge[](v12)(v11)
\Edge[](v13)(v14)
\Edge[](v13)(v7)
\Edge[](v14)(v7)
\Edge[](v1)(v13)
\Edge[](v3)(v14)
\end{tikzpicture}}}

%%%%%
% Figure C
\subfloat[A 7-vertex, $v$, with a 6-neighbor, $u_2$, and
two 5-neighbors, $u_1$ and $u_3$, with all pairs of $u_i$s non-adjacent.]{\makebox[.5\textwidth]{
\begin{tikzpicture}[scale = 10]
\tikzstyle{VertexStyle}=[shape = circle, minimum size = 6pt, inner sep = 1.2pt, draw]
\Vertex[x = 0.55, y = 0.55, L = \tiny {$v$}]{v0}
\Vertex[x = 0.40, y = 0.65, L = \tiny {$3$}]{v1}
\Vertex[x = 0.35, y = 0.50, L = \tiny {$u_1$}]{v2}
\Vertex[x = 0.70, y = 0.65, L = \tiny {$3$}]{v3}
\Vertex[x = 0.45, y = 0.40, L = \tiny {$1$}]{v4}
\Vertex[x = 0.65, y = 0.40, L = \tiny {$3$}]{v5}
\Vertex[x = 0.75, y = 0.50, L = \tiny {$u_3$}]{v6}
\Vertex[x = 0.55, y = 0.75, L = \tiny {$u_2$}]{v7}
\Vertex[x = 0.25, y = 0.60, L = \tiny {$1$}]{v8}
\Vertex[x = 0.25, y = 0.40, L = \tiny {}]{v9}
\Vertex[x = 0.40, y = 0.85, L = \tiny {$2$}]{v11}
\Vertex[x = 0.70, y = 0.85, L = \tiny {$2$}]{v12}
\Vertex[x = 0.85, y = 0.40, L = \tiny {}]{v13}
\Vertex[x = 0.85, y = 0.60, L = \tiny {}]{v14}
%\Vertex[x = 0.95, y = 0.50, L = \tiny {$w_3$}]{v15}
\Vertex[x = 0.55, y = 0.95, L = \tiny {}]{v16}
\Edge[](v0)(v1)
\Edge[](v0)(v2)
\Edge[](v0)(v4)
\Edge[](v0)(v5)
\Edge[](v0)(v6)
\Edge[](v0)(v3)
\Edge[](v3)(v6)
\Edge[](v6)(v5)
\Edge[](v2)(v1)
\Edge[](v2)(v4)
\Edge[](v4)(v5)
\Edge[](v7)(v0)
\Edge[](v7)(v3)
\Edge[](v7)(v1)
\Edge[](v9)(v2)
\Edge[](v9)(v4)
\Edge[](v8)(v2)
\Edge[](v8)(v1)
\Edge[](v8)(v9)
\Edge[](v12)(v3)
\Edge[](v12)(v7)
\Edge[](v11)(v1)
\Edge[](v11)(v7)
\Edge[](v14)(v6)
%\Edge[](v14)(v15)
\Edge[](v14)(v13)
\Edge[](v13)(v5)
\Edge[](v13)(v6)
%\Edge[](v13)(v15)
%\Edge[](v15)(v6)
\Edge[](v14)(v3)
\Edge[](v11)(v16)
\Edge[](v16)(v7)
\Edge[](v16)(v12)
\end{tikzpicture}}}
% Figure D
\subfloat[A 7-vertex, $v$, with a 5-neighbor, $u_1$, and
two 6-neighbors, $u_2$ and $u_3$, with all pairs of $u_i$s non-adjacent.]{\makebox[.5\textwidth]{
\begin{tikzpicture}[xscale = 10,yscale=10]
\tikzstyle{VertexStyle}=[shape = circle, minimum size = 6pt, inner sep = 1.2pt, draw]
\Vertex[x = 0.55, y = 0.55, L = \tiny {$v$}]{v0}
\Vertex[x = 0.40, y = 0.65, L = \tiny {$3$}]{v1}
\Vertex[x = 0.35, y = 0.50, L = \tiny {$u_1$}]{v2}
\Vertex[x = 0.70, y = 0.65, L = \tiny {$3$}]{v3}
\Vertex[x = 0.45, y = 0.40, L = \tiny {$1$}]{v4}
\Vertex[x = 0.65, y = 0.40, L = \tiny {$3$}]{v5}
\Vertex[x = 0.75, y = 0.50, L = \tiny {$u_3$}]{v6}
\Vertex[x = 0.55, y = 0.75, L = \tiny {$u_2$}]{v7}
\Vertex[x = 0.25, y = 0.60, L = \tiny {$1$}]{v8}
\Vertex[x = 0.25, y = 0.40, L = \tiny {}]{v9}
\Vertex[x = 0.40, y = 0.85, L = \tiny {$2$}]{v11}
\Vertex[x = 0.70, y = 0.85, L = \tiny {$2$}]{v12}
\Vertex[x = 0.85, y = 0.40, L = \tiny {$4$}]{v13}
\Vertex[x = 0.85, y = 0.60, L = \tiny {$4$}]{v14}
\Vertex[x = 0.95, y = 0.50, L = \tiny {$w_3$}]{v15}
\Vertex[x = 0.55, y = 0.95, L = \tiny {}]{v16}
\Edge[](v0)(v1)
\Edge[](v0)(v2)
\Edge[](v0)(v4)
\Edge[](v0)(v5)
\Edge[](v0)(v6)
\Edge[](v0)(v3)
\Edge[](v3)(v6)
\Edge[](v6)(v5)
\Edge[](v2)(v1)
\Edge[](v2)(v4)
\Edge[](v4)(v5)
\Edge[](v7)(v0)
\Edge[](v7)(v3)
\Edge[](v7)(v1)
\Edge[](v9)(v2)
\Edge[](v9)(v4)
\Edge[](v8)(v2)
\Edge[](v8)(v1)
\Edge[](v8)(v9)
\Edge[](v12)(v3)
\Edge[](v12)(v7)
\Edge[](v11)(v1)
\Edge[](v11)(v7)
\Edge[](v14)(v6)
\Edge[](v14)(v15)
\Edge[](v13)(v5)
\Edge[](v13)(v6)
\Edge[](v13)(v15)
\Edge[](v15)(v6)
\Edge[](v14)(v3)
\Edge[](v11)(v16)
\Edge[](v16)(v7)
\Edge[](v16)(v12)
\end{tikzpicture}}}

% Figure E
\subfloat[A 7-vertex, $v$, with a 5-neighbor, $u_2$, and
two 6-neighbors, $u_1$ and $u_3$, with all pairs of $u_i$s non-adjacent.]{\makebox[.5\textwidth]{
\begin{tikzpicture}[scale = 10]
\tikzstyle{VertexStyle}=[shape = circle, minimum size = 6pt, inner sep = 1.2pt, draw]
\Vertex[x = 0.55, y = 0.55, L = \tiny {$v$}]{v0}
\Vertex[x = 0.40, y = 0.65, L = \tiny {$2$}]{v1}
\Vertex[x = 0.35, y = 0.50, L = \tiny {$u_1$}]{v2}
\Vertex[x = 0.70, y = 0.65, L = \tiny {$3$}]{v3}
\Vertex[x = 0.45, y = 0.40, L = \tiny {$2$}]{v4}
\Vertex[x = 0.65, y = 0.40, L = \tiny {$3$}]{v5}
\Vertex[x = 0.75, y = 0.50, L = \tiny {$u_3$}]{v6}
\Vertex[x = 0.55, y = 0.75, L = \tiny {$u_2$}]{v7}
\Vertex[x = 0.25, y = 0.60, L = \tiny {$1$}]{v8}
\Vertex[x = 0.25, y = 0.40, L = \tiny {$1$}]{v9}
\Vertex[x = 0.15, y = 0.50, L = \tiny {$w_1$}]{v10}
\Vertex[x = 0.40, y = 0.85, L = \tiny {}]{v11}
\Vertex[x = 0.70, y = 0.85, L = \tiny {$2$}]{v12}
\Vertex[x = 0.85, y = 0.40, L = \tiny {$4$}]{v13}
\Vertex[x = 0.85, y = 0.60, L = \tiny {$4$}]{v14}
\Vertex[x = 0.95, y = 0.5, L = \tiny {$w_3$}]{v15}
\Edge[](v0)(v1)
\Edge[](v0)(v2)
\Edge[](v0)(v4)
\Edge[](v0)(v5)
\Edge[](v0)(v6)
\Edge[](v0)(v3)
\Edge[](v3)(v6)
\Edge[](v6)(v5)
\Edge[](v2)(v1)
\Edge[](v2)(v4)
\Edge[](v4)(v5)
\Edge[](v7)(v0)
\Edge[](v7)(v3)
\Edge[](v7)(v1)
\Edge[](v9)(v2)
\Edge[](v9)(v4)
\Edge[](v8)(v2)
\Edge[](v8)(v1)
\Edge[](v10)(v8)
\Edge[](v10)(v2)
\Edge[](v10)(v9)
\Edge[](v12)(v3)
\Edge[](v12)(v7)
\Edge[](v11)(v1)
\Edge[](v11)(v7)
\Edge[](v12)(v11)
\Edge[](v14)(v6)
\Edge[](v14)(v15)
\Edge[](v13)(v5)
\Edge[](v13)(v6)
\Edge[](v13)(v15)
\Edge[](v15)(v6)
\Edge[](v14)(v3)
\end{tikzpicture}}}

\caption{The five cases of Lemma~\ref{7566}.}
\label{fig7566}
\end{figure}
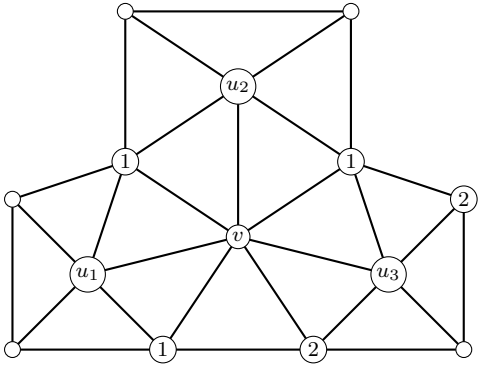
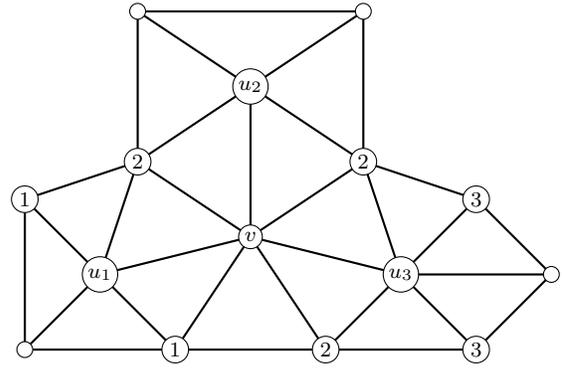
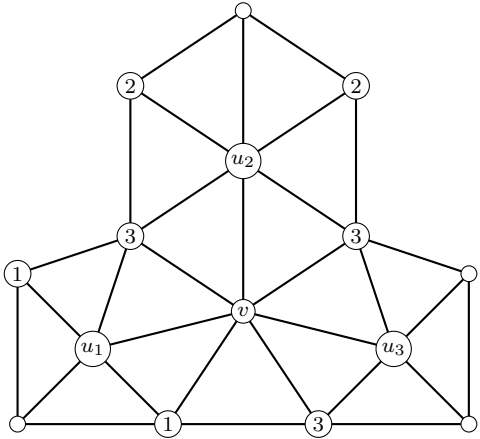
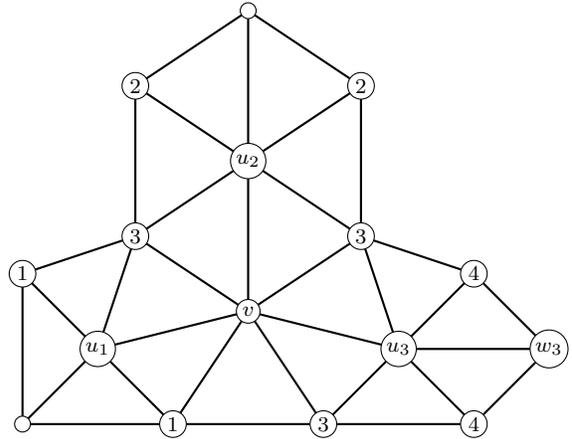
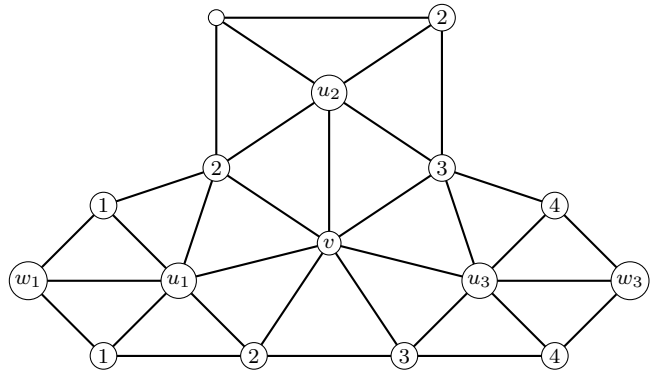
Figure~\ref{fig7566}(A) shows a $7$-vertex with three pairwise non-adjacent
$5$-neighbors.  Here, all pairs of vertices with numeric labels are at 
distance at most three, so they must be distinct.

In Figure~\ref{fig7566}(B), all pairs of vertices with numeric labels are again
at distance at most three, except for one vertex labeled $1$ which
is drawn at distance four from each vertex labeled $3$.  The only possible problem
is if one pair of vertices labeled $1$ and $3$ are actually the same vertex, while another pair
labeled $1$ and $3$ are adjacent; these pairs must be disjoint, since otherwise we
have a separating 3-cycle.  The pair that are adjacent must be drawn at
distance at least three, to avoid a separating 3-cycle.  Hence, we need only
consider the case where the vertices labeled $1$ and $3$ drawn at distance three are
adjacent, and the other pair labeled $1$ and $3$ are the same vertex $x$.  However, this is
impossible, since then the adjacent pair are on opposite sides of the cycle $u_1vu_3x$.

In Figure~\ref{fig7566}(C), all pairs of vertices with numeric labels are at
distance at most three, so they must be distinct.

Consider Figure~\ref{fig7566}(D).  None of the vertices labeled $3$ can
be the same as any other numerically labeled vertices since they are all distance at
most $3$ apart.  Similarly, none of the vertices labeled $1$ and $2$ can be the same.
So we need only consider the case that vertices labeled $4$
are the same as those labeled $1$ or $2$.  If a pair of vertices labeled $2$ and $4$ are
the same, then it must be the pair that are drawn at distance $4$; call this
vertex $x$.  In this case, we unlabel the vertices labeled $4$ and label $w_3$ with $3$.  Now, thanks
to $x$, all labeled vertices are distance at most three apart; hence, they
must be distinct.  So we may assume that vertices labeled $4$ are not the same as
those labeled $2$.  

Suppose instead that a vertex labeled $4$ is the same as
one labeled $1$; call this vertex $y$.  This is only a problem if also some
pair of vertices labeled 1 and 4 are adjacent.  But this is impossible as follows.
Since the pair of vertices labeled 1 have a common neighbor, they cannot be
adjacent; similarly for the pair labeled 4.
So the pairs that are identified and adjacent must be disjoint.  Further, the
identified pair must contain the rightmost vertex labeled 1.  If it is
identified with the bottom vertex labeled 4, then the remaining vertices cannot be
adjacent, since they are on opposite sides of the 4-cycle $u_1vu_3y$.
If it is identified with the top vertex labeled 4, then the remaining pair
cannot be adjacent, since they have a common neighbor.

Finally, consider Figure~\ref{fig7566}(E).  By horizontal symmetry (and planarity), we
assume that the vertices labeled $2$ that are drawn at distance three are indeed
non-adjacent; furthermore, we can assume that the vertices labeled $1$ and $2$
that are drawn at distance four are distinct.  If not, then we reflect across
the edge $u_2v$.  Hence, in forming $G'$ we can contract the vertices labeled
$2$ to a single vertex (we can also contract the vertices labeled $3$ to a
single vertex).  So we only need to consider the vertices labeled $1$ and $4$. 
The only possible problem is if some pair of vertices labeled $1$ and $4$
that are drawn at distance four are actually the same vertex $x$.  Suppose
this is the case.  If $w_1$ and $w_3$ are distinct, then we
neglect the vertices labeled $1$ and $4$ altogether;
instead we label $w_1$ as $2$ and $w_3$ as $3$.  Due to $x$,
all pairs of vertices with numeric labels are now distance at most three.
(Also, we can assume that $w_1$ is not adjacent to the vertex labeled 2 that is
drawn at distance 4; if not, then we again reflect across edge $u_1v$.)
So assume that $w_1$ and $w_3$ are identified.  %This means that the top
Now we switch the vertex identifications we use to form $G'$.
Contract the two vertices labeled 4 onto $u_3$; also contract onto $u_1$ its two
neighbors labeled 2, the topmost vertex labeled 3, and $w_{1/3}$.  As usual,
we color this smaller graph by minimality; when we lift this coloring to $G$, vertex
$v$ and each vertex $u_i$ has enough available colors that we can extend the
coloring by Fact~1.  This finishes Figure~\ref{fig7566}(E) and completes the
proof of the lemma.  
\end{proof}

Now we use discharging to prove that every planar graph has a 2-fold 9-coloring.

\begin{maintheorem}
Every planar graph $G$ has a 2-fold 9-coloring.  In particular,
$\chif(G)\le\frac92$.
\end{maintheorem}

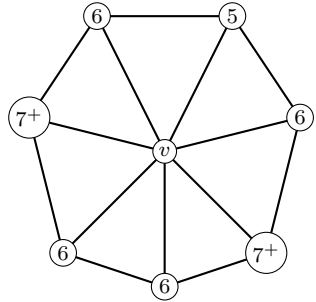
\begin{figure}[h!]
\centering
\begin{tikzpicture}[scale = 9]
\tikzstyle{VertexStyle}=[shape = circle, minimum size = 6pt, inner sep = 1.2pt, draw]
\Vertex[x = 0.60, y = 0.70, L = \tiny {$v$}]{v0}
\Vertex[x = 0.50, y = 0.90, L = \tiny {$6$}]{v1}
\Vertex[x = 0.70, y = 0.90, L = \tiny {$5$}]{v2}
\Vertex[x = 0.80, y = 0.75, L = \tiny {$6$}]{v3}
\Vertex[x = 0.40, y = 0.75, L = \tiny {$7^+$}]{v4}
\Vertex[x = 0.45, y = 0.55, L = \tiny {$6$}]{v5}
\Vertex[x = 0.60, y = 0.50, L = \tiny {$6$}]{v6}
\Vertex[x = 0.75, y = 0.55, L = \tiny {$7^+$}]{v7}
\Edge[](v1)(v0)
\Edge[](v2)(v0)
\Edge[](v3)(v0)
\Edge[](v4)(v0)
\Edge[](v5)(v0)
\Edge[](v6)(v0)
\Edge[](v7)(v0)
\Edge[](v5)(v6)
\Edge[](v5)(v4)
\Edge[](v1)(v4)
\Edge[](v1)(v2)
\Edge[](v2)(v3)
\Edge[](v3)(v7)
\Edge[](v7)(v6)
\end{tikzpicture}
\caption{A 7-vertex $v$ gives no charge to any crowded 5-neighbor.}
\label{fig:whycrowded}
\end{figure}

\begin{proof}
The second statement follows from the first, which we prove now.
Let $G$ be a minimal counterexample to the theorem.
We will use the discharging method with initial charge $d(v)-6$ for each vertex $v$.  
We write $\ch(v)$ to denote the initial charge and $\ch^*(v)$
to denote the charge after redistributing.
By Euler's Formula, $\sum_{v\in V(G)} \ch(v) = -12$.  By assuming that $G$
satisfies the conditions stipulated in Lemmas~\ref{noTriangles}--\ref{7566}, we
redistribute the charge (without changing its sum) so that every
vertex finishes with nonnegative charge. 
This yields the obvious contradiction $-12 = \sum_{v\in V(G)}\ch(v) =
\sum_{v\in V(G)}\ch^*(v) \ge 0$.

We need a few definitions.
For a vertex $v$, let $H_v$ denote the subgraph induced by the $5$-neighbors and
$6$-neighbors of $v$.  If some $w\in V(H_v)$ has $d_{H_v}(w)=0$, then $w$
is an \emph{isolated} neighbor of $v$; otherwise $w$ is a \emph{non-isolated}
neighbor.  A non-isolated 5-neighbor of a vertex $v$ is \emph{crowded} (with
respect to $v$) if it has two $6$-neighbors in $H_v$.  We use crowded 5-neighbors
in the discharging proof to help ensure that $7$-vertices finish with sufficient
charge, specifically to handle the configuration in Figure~\ref{fig:whycrowded}.
We redistribute charge via the following four rules; they are applied
simultaneously, wherever applicable.

\begin{enumerate}
\item[(R1)] Each $8^+$-vertex gives charge $\frac12$ to each isolated
$5$-neighbor and charge $\frac14$ to each non-isolated $5$-neighbor.
\item[(R2)] Each $7$-vertex gives charge $\frac12$ to each isolated
$5$-neighbor, charge $0$ to each crowded $5$-neighbor and charge $\frac14$ to
each remaining $5$-neighbor.
\item[(R3)] Each $7^+$-vertex gives charge $\frac14$ to each $6$-neighbor.
\item[(R4)] Each $6$-vertex gives charge $\frac12$ to each $5$-neighbor.
\end{enumerate}
\bigskip

To show that every vertex $v$ finishes with nonnegative charge,
we consider $d(v)$.
\bigskip

$\bf d(v)\ge 8$: 
We will show that $v$ gives away charge at most $\frac{d(v)}4$.
Since $d(v) \ge 8$, we have $\ch(v) = d(v)-6\ge \frac{d(v)}4$, so this will
imply $\ch^*(v)\ge 0$.
Rather than giving away charge by rules (R1) and (R3), instead let $v$ give
charge $\frac14$ to each neighbor.  Now let each isolated $5$-neighbor $w$ take
also the charge $\frac14$ that $v$ gave to the neighbor that clockwise around
$v$ succeeds $w$.  Now each neighbor of $v$ has received at least as much
charge as by rules (R1) and (R3) and $v$ has given away
charge $\frac{d(v)}{4}$.   Thus, when $v$ gives away charge according to rules
(R1) and (R3), this charge is at most $\frac{d(v)}4$, so $\ch^*(v)\ge 0$.

$\bf d(v)=7$:
First, suppose that $v$ has an isolated $5$-neighbor $w$. Let $x,y \in N(v)$ be
the two $7^+$-vertices that are common neighbors of $v$ and $w$. 
We will show that the total charge that $v$ gives to $N(v)\setminus \{x,y\}$ is at most
$\frac12$.  By Lemma~\ref{7566}, these four remaining vertices include at most
two $6^-$-vertices.  So, if $v$ gives them a total of more than $\frac12$, then
one of them must be another isolated $5$-neighbor.
But now the final $6^-$-vertex must be at distance $2$ from each of the
previous 5-neighbors, violating Lemma~\ref{7566}.
 
So instead assume that $v$ has no isolated $5$-neighbors.  Thus, if $v$ loses total charge
more than 1, then it must have at least five $6^-$-neighbors that receive charge
from it (since they each take charge $\frac14$).  So assume that $|H_v|\ge 5$. 
This implies that $H_v$ consists of either (i) a 7-cycle or (ii) a single path
or (iii) two paths.  Recall from Lemma~\ref{666}, that no 6-vertex has
non-adjacent $6^-$-neighbors.  This means that every vertex of degree 2 in $H_v$
is a 5-vertex; in other words, every vertex on a cycle or in the interior of a
path in $H_v$ is a 5-vertex.  

Now in each of cases (i)--(iii), $H_v$ has an independent 3-set containing at
least one 5-vertex; the only exception is if $H_v$ consists of a path on two
vertices and a path on three vertices, and the only 5-vertex is the internal
vertex on the longer path.  However, in this case the 5-vertex is a crowded
neighbor of $v$, as in Figure~\ref{fig:whycrowded},  so it receives no charge
from $v$.  Thus, $\ch^*(v)\ge 0$.

$\bf d(v)=6$:  
By Lemma~\ref{666}, we know that $v$ has at most two $6^-$-neighbors (and if
exactly two, then they are adjacent).  Now (R3) implies that $\ch^*(v)\ge
0+4(\frac14)-2(\frac12)=0$.

$\bf d(v)=5$:  
If $v$ has at least two $6$-neighbors, then $\ch^*(v)\ge -1+2(\frac12)=0$; so assume
that $v$ has at most one $6$-neighbor.
Now if $v$ has at least four $6^+$-neighbors, then 
%$v$ finishes with charge at least 
$\ch^*(v)\ge-1+4(\frac14)=0$ (since $v$ has at most one $6$-neighbor, $v$ is not
a crowded neighbor for any of its $7$-neighbors); so $v$ must have at least
two $5$-neighbors.
By Lemma~\ref{556}, these $5$-neighbors must be adjacent and $v$ has no
$6$-neighbors.  But now one of $v$'s three $7^+$-neighbors sees $v$ as an isolated
$5$-neighbor, so sends $v$ charge $\frac12$.  Thus, $\ch^*(v)\ge -1+\frac12+2(\frac14)=0$.  
%finishes with charge at least
This completes the proof.
\end{proof}

A natural question is whether our theorem could be strengthened to show that
every planar graph has a $t$-fold $s$-coloring, for some pair $(s,t)$ with
$\frac{s}{t}<\frac92$.  Clearly, such results are true for every pair $(s,t)$
with $\frac{s}t\ge 4$, since they follow from the 4 Color Theorem (this is
immediate since the Kneser graph $K_{s:t}$ contains $K_4$).
However,
here we note that a proof of any such result must differ significantly from the
proof of the Main Theorem.  In particular, we show that none of our reducibility
proofs, with the exceptions of those for separating triangles and vertices of
degree 4, remain valid for any pair $(s,t)$ with $\frac{s}{t}<\frac92$.
Recall that the proofs of Lemmas~\ref{556}--\ref{7566} all crucially relied on
Fact~\ref{fact1}.  Here we show that to prove an analogue of this fact, even
for $K_{1,2}$ (rather than $K_{1,3}$) requires that $\frac{s}t\ge \frac 92$.

Suppose that a copy of $K_{1,2}$ has a $t$-fold coloring whenever the leaves,
$u_1$ and $u_2$, are given lists of size $b$ and the center vertex, $v$, is
given a list of size $a$.
Consider the list assignment $L(u_1)=\{1,\ldots,b\}$,
$L(u_2)=\{a-b+1,\ldots,a\}$, and $L(v)=\{1,\ldots,a\}$.  Every $t$-fold
coloring from these lists uses at most $b-(a-b)=2b-a$ common colors on $u_1$ and
$u_2$, so uses at least $2t-(2b-a)$ distinct colors on $u_1$ and $u_2$.
So, to color $v$, we must have $a -(2t-(2b-a))\ge t$, which means $b\ge
\frac32t$.  Now consider an analogue of Lemma~\ref{556}, \ref{666}, or
\ref{7566} for $t$-fold $s$-coloring.  First we contract, color the smaller
graph by minimality, and lift the coloring to $G$.  Now $v$ has list size
$s-2t$, and each $u_i$ has list size $s-3t$.  So we have $a=s-2t$ and $b=s-3t$.
Now $s-3t=b\ge \frac32t$, so $\frac{s}t\ge \frac92$.

\section*{Acknowledgments}
As we mentioned in the introduction, the ideas in this paper come largely from
Albertson's proof~\cite{Albertson} that planar graphs have independence ratio at least $\frac29$.
In fact, many of the reducible configurations that we use here are special cases
of the reducible configurations in that proof.  We very much like that
paper, and so it was a pleasure to be able to extend Albertson's work.
%It seems that the part of his own proof that Albertson was least pleased with
%was verifying ``unavoidability'', i.e., showing that every planar graph contains a
%reducible configuration.  In the introduction to~\cite{Albertson}, he wrote:
%``Finally Section 4 is devoted to a massive, ugly edge counting which
%demonstrates that every planar triangulation of the plane must contain some
%forbidden subgraph.''  In the appendix that follows, we give a short proof of
%this same unavoidability statement, via discharging.  We think Mike might
%have liked it.

The first author thanks his Lord and Savior, Jesus Christ.

\bibliographystyle{abbrvplain}
\bibliography{45ct}

\begin{thebibliography}{10}

\bibitem{Albertson}
M.~O. Albertson.
\newblock A lower bound for the independence number of a planar graph.
\newblock {\em J. Combinatorial Theory Ser. B}, 20(1):84--93, 1976.

\bibitem{berge}
C.~Berge.
\newblock {\em Graphes et hypergraphes}.
\newblock Dunod, Paris, 1970.
\newblock Monographies Universitaires de Math{\'e}matiques, No. 37.

\bibitem{birkhoff}
G.~Birkhoff.
\newblock The reducibility of maps.
\newblock {\em American Journal of Mathematics}, pages 115--128, 1913.

\bibitem{Diestel}
R.~Diestel.
\newblock {\em Graph theory}, volume 173 of {\em Graduate Texts in
  Mathematics}.
\newblock Springer-Verlag, Berlin, third edition, 2005.

\bibitem{Hell79}
P.~Hell.
\newblock An introduction to the category of graphs.
\newblock In {\em Topics in graph theory ({N}ew {Y}ork, 1977)}, volume 328 of
  {\em Ann. New York Acad. Sci.}, pages 120--136. New York Acad. Sci., New
  York, 1979.

\bibitem{HRS}
A.~J.~W. Hilton, R.~Rado, and S.~H. Scott.
\newblock A ({$<5$})-colour theorem for planar graphs.
\newblock {\em Bull. London Math. Soc.}, 5:302--306, 1973.

\bibitem{Kainen}
P.~C. Kainen.
\newblock A generalization of the {$5$}-color theorem.
\newblock {\em Proc. Amer. Math. Soc.}, 45:450--453, 1974.

\bibitem{Lovasz78}
L.~Lov{\'a}sz.
\newblock Kneser's conjecture, chromatic number, and homotopy.
\newblock {\em J. Combin. Theory Ser. A}, 25(3):319--324, 1978.

\bibitem{Naserasr06}
R.~Naserasr.
\newblock {$K_5$}-free bound for the class of planar graphs.
\newblock {\em European J. Combin.}, 27(7):1155--1158, 2006.

\bibitem{NOdM06}
J.~Ne{\v{s}}et{\v{r}}il and P.~Ossona~de Mendez.
\newblock Folding.
\newblock {\em J. Combin. Theory Ser. B}, 96(5):730--739, 2006.

\bibitem{Sauer01}
N.~Sauer.
\newblock {Hedetniemi's Conjecture---A Survey}.
\newblock {\em Discrete Math.}, 229(1-3):261--292, 2001.
\newblock Combinatorics, graph theory, algorithms and applications.

\bibitem{SU-book}
E.~R. Scheinerman and D.~H. Ullman.
\newblock {\em {Fractional Graph Theory}}.
\newblock Wiley-Interscience Series in Discrete Mathematics and Optimization.
  John Wiley \& Sons, Inc., New York, 1997.
\newblock A rational approach to the theory of graphs, With a foreword by
  Claude Berge, A Wiley-Interscience Publication.

\bibitem{Wagner37}
K.~Wagner.
\newblock \"{U}ber eine {E}igenschaft der ebenen {K}omplexe.
\newblock {\em Math. Ann.}, 114(1):570--590, 1937.

\end{thebibliography}
\end{document}